\documentclass[12pt]{amsart}
\usepackage{amscd,amssymb}
\usepackage[graph,frame,poly,arc]{xy}  
\usepackage[plainpages,backref,urlcolor=blue]{hyperref}

\theoremstyle{plain}
\newtheorem{thm}[subsection]{Theorem}

\newtheorem{prop}[subsection]{Proposition}
\newtheorem{cor}[subsection]{Corollary}

\theoremstyle{definition}
\newtheorem{rk}[subsection]{Remark}
\newtheorem{definition}[subsection]{Definition}
\newtheorem{ex}[subsection]{Example}
\newtheorem{conj}[subsection]{Conjecture}
\newtheorem{question}[subsection]{Question}

\numberwithin{equation}{section}
\newcommand{\A}{{\mathcal A}}
\newcommand{\B}{{\mathcal B}}
\newcommand{\G}{{\mathcal G}}

\newcommand{\M}{{\mathcal M}}

\newcommand{\CC}{{\mathcal C}}
\newcommand{\LL}{{\mathcal L}}
\newcommand{\E}{{\mathcal E}}

\newcommand{\Q}{\mathbb{Q}}

\newcommand{\C}{\mathbb{C}}
\newcommand{\PP}{\mathbb{P}}

\newcommand{\N}{\mathbb{N}}

\DeclareMathOperator{\codim}{codim}

\DeclareMathOperator{\depth}{depth }

\begin{document}

\title [Numerical invariants and moduli spaces for line arrangements]
{Numerical invariants and moduli spaces for line arrangements}

\author[A. Dimca]{Alexandru Dimca$^1$}
\address{Universit\'e C\^ ote d'Azur, CNRS, LJAD, France  }
\email{dimca@unice.fr}

\author[D. Ibadula]{Denis~Ibadula}
\address{Ovidius University, Faculty of Mathematics and Informatics, 124 Mamaia Blvd., 900527 Constan\c{t}a, Romania}
\email{denis.ibadula@univ-ovidius.ro}

\author[A. D. M\u acinic]{Daniela Anca~M\u acinic}
\address{Simion Stoilow Institute of Mathematics, 
P.O. Box 1-764,
RO-014700 Bucharest, Romania}
\email{Anca.Macinic@imar.ro}

\thanks{$^1$ Partially supported by Institut Universitaire de France.}



\begin{abstract} 
Using several numerical invariants, we study a partition of the space of line arrangements in the complex projective plane, given by the intersection lattice types. We offer also a new characterization of the free plane curves using the Castelnuovo-Mumford regularity of the associated Milnor/Jacobian algebra.

\end{abstract}
 
\maketitle


\section{Introduction} 
Line arrangements in the complex projective plane $\PP^2$ look like being simple objects, but a lot of questions related to them are still unanswered, e.g. Terao's conjecture saying that the freeness of such an arrangement is determined by the combinatorics, see Conjecture \ref{conjT} below for a statement, \cite[Chapter 8]{DHA} for more information, and \cite{Yo} for a survey. Or the conjecture that the monodromy of the associated Milnor fiber is determined by the combinatorics, see \cite{PS} for a recent result and \cite{DHA} for more information.

In order to treat such questions, the study of parameter spaces (a.k.a moduli spaces) of line arrangements has being developed, centered especially on the irreducibility/connectivity questions, see  \cite{10lines},  \cite{10lines'}, \cite{Artal06}, \cite{NY},  \cite{Ye}.

In this paper, the new idea is to look at the way in which the parameter spaces $A(L)$ and $X(L)$
of line arrangements with a given intersection lattice $L$ behave when the lattice $L$ changes. In section 2 we describe two parameter spaces for the line arrangements $\A:f=0$ in $\PP^2$ having $d$ lines, namely $A(d)$ and $X(d)$, which are both smooth irreducible varieties, see Corollary \ref{corAnew}. To partition these two varieties $A(d)$ and $X(d)$ into finer strata, keeping track of the properties of the line arrangements, we use several numerical invariants and study their semi-continuity properties in Proposition \ref{prop2}. We consider in this section both line arrangements and arbitrary reduced curves in $\PP^2$, in order to point out that the numerical invariants associated to line arrangements enjoy special properties, see Corollary \ref{corArr}.

In section 3 we recall the definition and main properties of free and nearly free plane curves. Then we prove that a classical invariant in Commutative Algebra, namely the Castelnuovo-Mumford regularity, coincides, when applied to the Milnor/Jacobian  algebra $M(f)$ of a reduced plane curve $C:f=0$, to a naive invariant $st(f)$, coming from the Hilbert function of the graded algebra $M(f)$, exactly when the curve $C$ is free, see Theorem \ref{thmA} and Corollaries \ref{corG} and \ref{corB}. Corollary \ref{corB} depends on a key result due to H. Schenck, see \cite{SchCMH}. We end this section by noting that our partitions of the spaces $A(d)$ and $X(d)$ are $G-$equivariant, where $G=Aut(\PP^2)$ is acting in the obvious way on these parameter spaces. In Proposition \ref{prop3} we describe the dimension of a line arrangement under this $G-$action. 

In section 4 we fix an integer $d \geq 4$ and denote by $\LL(d)$ the set of all possible intersection lattices of arrangements with $d$ lines, modulo lattice isomorphisms. For each
(isomorphism class of) lattice $L \in \LL(d)$, we denote by $X(L)$ the subset of $X(d)$ consisting of line arrangements having an intersection lattice isomorphic to $L$.
Hence the subsets $X(L)$ for $L \in \LL(d)$ form the strata of a partition of the smooth irreducible variety $X(d)$. The first properties of these strata are given in Proposition \ref{prop4}. Then we discuss several examples of simple lattices $L$ and of corresponding strata $X(L)$, e.g. the lattice $L_{gen}$ corresponding to the generic line arrangement is discussed in Example \ref{ex1} and an obvious generalization, the lattice $L(d,m)$, is considered in Proposition \ref{prop5}. Other lattices occurs in Proposition \ref{propAD}, Example \ref{exrkAD}, Example \ref{exrkexpo}, the last two describing free (resp. nearly free) line arrangements. Note that our results on a stratum $X(L)$ (e.g. dimension, connectivity) easily translate into properties of the quotient $X(L)/G$. We also reprove in a new way a result due to Toh\u aneanu \cite{T}, giving the classification of line arrangements with a Jacobian syzygy of minimal degree 2, see Theorem \ref{thmToh}.

In section 5 we point out the complexity of the stratification $(X(L))_{L \in \LL(d)}$ of the space $X(d)$. First we describe all the strata when $d=4,5,6$, and pay particular attention to the strata formed by (nearly) free arrangements. We explain just after Proposition \ref{propd=6} that these stratification do not satisfy the frontier condition in general, in particular they are not Whitney regular stratifications. In Remark \ref{rkd=7} and in the answer to Question \ref{qd=7} we show that some nice features of this stratification noticed when $d \leq 6$ do not extend to higher degrees $d$.

In the final section we discuss Terao's conjecture in the case of line arrangements, recall the known results and give a new proof for Theorem \ref{thmT2}. Finally, in Proposition \ref{proptau}  we give a generalization of the result saying that the generic arrangement is not free when $d \geq 4$.

\bigskip

The authors would like to thank the Oberwolfach Research Institute for Mathematics, where the major part on the work on this project was done during a RIP program.

We also thank Torsten Hoge for his useful remarks on the previous version of this paper.

\section{General facts on plane curves and line arrangements} 

\subsection{Two parameter spaces for line arrangements: $A(d)$ and $X(d)$} 

Let $S=\C[x,y,z]$ be the graded polynomial ring in the variables $x,y,z$ with complex coefficients, and $S_m$ be the vector space of degree $m$ homogeneous polynomials in $S$. Fix an integer $d \geq 1$ and regard the projective space $\CC(d)=\PP(S_d)$ as the parameter space of degree $d$ curves in $\PP^2$. 

 \begin{prop}
\label{prop1} 
The subset $\CC(d)_0 \subset \CC(d)$ corresponding to curves having only isolated singularities is  Zariski open and dense in $\CC(d)$. The subset $A(d) \subset \CC(d)_0$ corresponding to  line arrangements $\A$ consisting of $d$ distinct lines is a closed Zariski subset in $\CC(d)_0$.
\end{prop}

\proof
For the first claim, note that the complement $\CC(d) \setminus \CC(d)_0$ is the union of the finite family of Zariski closed subsets given by the images of the obvious mappings
$$\phi_m: \PP(S_m) \times \PP(S_{d-2m}) \to \PP(S_d),$$
$(A,B) \mapsto A^2 B$, for $m=1,2,..., [d/2]$.
For the second claim, we consider the map
\begin{equation}
\label{keymap1}
\psi: \PP(S_1)^d \to \PP(S_d),
\end{equation}
given by $(\ell_1,...,\ell_d) \mapsto \ell_1 \cdot \ell_2 \cdot ... \cdot \ell_d$.
Then $A(d)$ is just the intersection of the set $\CC(d)_0$ with the image of the mapping $\psi$.
\endproof
Note that 
\begin{equation}
\label{keymap2}
X(d):=\psi^{-1}(A(d))
\end{equation}
is exactly the set of linear forms $(\ell_1,...,\ell_d) \in  \PP(S_1)^d$ such that $\ell_i \ne \ell_j$ for $i \ne j$, and the restriction $\psi: X(d) \to A(d)$ is a Galois covering with structure group the symmetric group $\sigma_d$ on $d$ elements.

\begin{cor}
\label{corAnew}
The parameter spaces $A(d)$ and $X(d)$ are smooth, irreducible algebraic varieties of dimension $2d$. The space $X(d)$ is simply connected and the fundamental group of the space $A(d)$ is given by
$$\pi_1(A(d))= \sigma_d.$$
\end{cor}

\proof
The only claim that needs some explanation is the fact that $X(d)$ is simply connected. This follows from the fact that $X(d)$ is obtained from the simply connected variety $\PP(S_1)^d$ by removing the codimension 2 linear subvarieties $\Delta_{ij}: \ell_i=\ell_j$ for all $1 \leq i<j \leq d$.
\endproof

\subsection{Some numerical invariants for plane curves and line arrangements}

For a polynomial $f \in S_d$, we denote by $J_f$ the ideal in $S$ generated by the partial derivatives $f_x,f_y,f_z$, and call $J_f$ the Jacobian ideal of $f$. The graded ring $M(f)=S/J_f$ is called the Milnor or Jacobian algebra of $f$.
 We define
\begin{equation}
\label{m_k}
m_k(f)=\dim M(f)_k.
\end{equation}
Note that one has $m_k(f)=\tau(f)$ for $k > 3(d-2)$ and $f \in \CC(d)_0$, where
$\tau(f)$ is the total Tjurina number of the reduced plane curve $C(f): f=0$, see \cite{CD}. We also denote by $\tau(C(f))$ this number, and note that it is nothing else but the degree of the Jacobian ideal $J_f$.

The minimal degree of a Jacobian syzygy for $f$ is the integer $mdr(f)$
defined to be the smallest integer $r \geq 0$ such that there is a nontrivial relation
\begin{equation}
\label{rel_m}
 af_x+bf_y+cf_z=0
\end{equation}
among the partial derivatives $f_x, f_y$ and $f_z$ of $f$ with coefficients $a,b,c$ in $S_r$.
We denote by $AR(f)$ the graded $S$-module consisting of all the triples $(a,b,c) \in S^3$ satisfying  \eqref{rel_m}. In fact $AR(f)$ depends only on the class of $f\in S_d$ in $\CC(d)=\PP(S_d)$.

 \begin{prop}
\label{prop2} 

\begin{enumerate}

\item 
The subset $\{ f \in S_d \ : m_k(f) \leq m\}  \subset \CC(d)$  is  Zariski open and dense in $\CC(d)$ for any $k\geq 0$ and any $m \geq 0$. In particular, the following two sets 
$\{ f \in \CC(d)_0 \ : \tau(f) \leq m\}  \subset \CC(d)_0$  
and $\{ f \in A(d) \ : \tau(f) \leq m\}  \subset A(d)$,
are  Zariski open and dense in $\CC(d)_0$ (resp. in $A(d)$) for any $k\geq 0$ and any $m \geq 0$. 

\item The subset $\{ f \in S_d \ : mdr(f) \leq m\}  \subset \CC(d)$  is  Zariski closed in $\CC(d)$ for any $m\geq 0$. In particular, the following two sets 
$\{ f \in \CC(d)_0 \ : mdr(f) \leq m\}  \subset \CC(d)_0$  
and $\{ f \in A(d) \ : mdr(f) \leq m\}  \subset A(d)$,
are  Zariski closed in $\CC(d)_0$ (resp. in $A(d)$) for any $m \geq 0$.

\end{enumerate}

\end{prop}

\proof The first claim is clear by the semicontinuity properties of the rank of a matrix,
whose rows are obtained by taking all the coefficients of the polynomials $\mu g$,
where $\mu$ runs through the set of monomials of degree $k-d+1$ in $S$ and $g \in \{f_x,f_y,f_z\}$. 

To prove the second claim, consider the closed subvariety $Y_m$ in $\PP(S_m^3)\times \PP(S_d)$ given by
$$Y_m=\{((a,b,c),f) \  \ : \ \ af_x+bf_y+cf_z=0\}.$$
Note that a polynomial $f \in S_d$ satisfies $mdr(f) \leq m$ if and only if $[f] \in \PP(S_d)$ is in the image of $Y_m$ under the second projection. 

\endproof

\begin{definition}
\label{def}
For a polynomial $f \in \CC(d)_0$, we recall the following invariants.

\noindent (i) the {\it coincidence threshold} 
$$ct(f)=\max \{q:\dim M(f)_k=\dim M(f_s)_k \text{ for all } k \leq q\},$$
with $f_s$  a homogeneous polynomial in $S$ of the same degree $d$ as $f$ and such that $C_s:f_s=0$ is a smooth curve in $\PP^2$.

\noindent (ii) the {\it stability threshold} 
$st(f)=\min \{q~~:~~\dim M(f)_k=\tau(f) \text{ for all } k \geq q\}.$

\noindent (iii) the {\it regularity} 
$reg(f)$ is the Castelnuovo--Mumford regularity of the Milnor algebra $M(f)$, regarded as a graded $S$-module, see \cite[Chapter 4]{Eis}.
\end{definition}

The exact sequences 
\begin{equation} 
\label{EXSE}
 0 \to AR(f)(-(d-1)) \to S(-(d-1))^3 \xrightarrow{(f_x,f_y,f_z)} S \to M(f) \to 0
\end{equation} 
and
\begin{equation} 
\label{EXSEs}
 0 \to AR(f_s)(-(d-1)) \to S(-(d-1))^3 \xrightarrow{(f_{s,x},f_{s,y},f_{s,z})} S \to M(f_s) \to 0
\end{equation} 
and the fact that $AR(f_s)_k=0$ for $k<d-1$ imply that
\begin{equation} 
\label{REL}
ct(f) \geq mdr(f)+d-2,
\end{equation} 
with equality for $mdr(f) <d-1.$ To have equality always, it is convenient to introduce the invariant
$mdr_e(f)$, the minimal degree of an essential Jacobian relation for $f$, which is by definition the minimal degree of a relation \eqref{rel_m}, where the triple $(a,b,c)$ does not belong to the $S$-submodule of $AR(f)$ generated by the Koszul relations $(f_y,-f_x,0)$, $(f_z,0,-f_x)$ and $(0,f_z,-f_y)$. With this definition we always have
\begin{equation} 
\label{REL_e}
ct(f)= mdr_e(f)+d-2,
\end{equation} 
see \cite{DStEdin}.
 \begin{cor}
\label{corArr} With the above notation, the following hold.
\begin{enumerate}

\item Let $C:f=0$ be a singular, reduced plane curve of degree $d\geq 3 $ in $\PP^2$. Then
$$mdr(f) \leq d-1, \  mdr_e(f) \leq 2(d-2)\text{  and  }ct(f)\leq 3(d-2).$$
 Moreover, if $\tau(C)=1$, all
these inequalities are equalities.

\item Let $\A:f=0$ be an arrangement having $d\geq 2$ lines. Then
$$mdr(f)=mdr_e(f) \leq d-2 \text{ and } ct(f) \leq 2(d-2).$$
Moreover, both of these inequalities are equalities for a generic arrangement.
\end{enumerate}
\end{cor}
\proof
To prove (1), note that $mdr(f) \leq d-1$ follows from the existence of Koszul relations.
The inequality $mdr_e(f) \leq 2(d-2)$ follows using \eqref{REL_e} and the obvious fact that
$ct(f) \leq T=3(d-2).$ When $\tau(C)=1$, it follows from \cite[Example 4.3]{DStEdin}
that $ct(f)=T$. Moreover, such a curve is nodal and irreducible, and hence $mdr(f) \geq d-1$, by
\cite[Theorem 4.1]{DStEdin}.
This remark completes the proof of the first claim.

In view of Proposition \ref{prop2} (2), to prove (2) it is enough to show that $mdr(f) =d-2$
when $\A$ is a generic arrangement. But this follows from 
\cite[Theorem 4.1]{DStEdin}. 
\endproof

The invariants $st(f)$ and $reg(f)$ are closely related, as  the next section shows. However, they do not seem to satisfy semicontinuity properties similar to those in Proposition \ref{prop2}, see Remark \ref{rksemic}.

\section{Free and nearly free plane curves}

\subsection{Free plane curves}
For the equivalence of the properties in the next definition, we refer to \cite{ST}. See also \cite[Chapter 8]{DHA}.
\begin{definition} 
\label{deffree}
The curve $C:f=0$ is a { free divisor} if the following  equivalent conditions hold.

\begin{enumerate}

\item The Milnor algebra $M(f)$ is a Cohen-Macaulay $S$-module.

\item The minimal graded resolution of the Milnor algebra $M(f)$ as an $S$-module
has the following  form
$$0 \to S(-d_1-d+1) \oplus S(-d_2-d+1) \to S^3(-d+1) \xrightarrow{(f_x,f_y,f_z)}  S$$
for some positive integers $d_1, d_2$.
\item The graded $S$-module $AR(f)$ is free of rank 2, i.e. there is an isomorphism 
$$AR(f)=S(-d_1) \oplus S(-d_2)$$
for some positive integers $d_1, d_2$.
\end{enumerate}

\end{definition}
When $C$ is a free divisor, the integers $d_1 \leq d_2$ are called the {exponents} of $C$.  They satisfy the relations 
\begin{equation}
\label{free1}
 d_1+d_2=d-1 \text{ and } \tau(C)=(d-1)^2 - d_1d_2,
\end{equation}
where $\tau(C)$ is the total Tjurina number of $C$, see for instance \cite{DS14}. For a free curve, one  has $mdr(f)=d_1$, $ct(f)=d_1+d-2$ in view of \eqref{REL}, and $st(f)=d_2+d-3$, see for instance \cite{Dmax}.
\begin{definition} 
\label{def2}

The curve $C:f=0$ is a { nearly free divisor} 
if  the Milnor algebra $M(f)$ has a minimal graded resolution of the form
$$0 \to S(-d-d_2) \to S(-d-d_1+1) \oplus S^2(-d-d_2+1) \to S^3(-d+1) \xrightarrow{(f_x,f_y,f_z)}  S$$
for some integers $1 \leq d_1 \leq d_2$, called the exponents of $C$.

\end{definition}
For a nearly free curve, the exponents satisfy $d_1+d_2=d$, and one  has $mdr(f)=d_1$, $ct(f)=d_1+d-2$ in view of \eqref{REL}, and $st(f)=d_2+d-2$ by the results in \cite{Dmax}. 

 \begin{thm}
\label{thmA} Let $C:f=0$ be a reduced plane curve.
 Then 
 $$st(f)-1 \leq reg(f) \leq st(f),$$
 and the equality  $reg(f)=st(f)$ holds if and only if $C:f=0$ is a free curve.

\end{thm}

\proof

 Let $H_{M(f)}$ (resp. $P_{M(f)}$) be the Hilbert function (resp. the Hilbert polynomial) of the graded $S$-module $M(f)$. Then \cite[Theorem 4.2]{Eis} implies that
$$H_{M(f)}(k)= P_{M(f)}(k)$$
for any $k \geq reg(f)+1$. Since for a reduced plane curve the Hilbert polynomial $P_{M(f)}$ is just the constant $\tau(C)$, it follows from the definition of the stability threshold $st(f)$ that $st(f) \leq reg(f)+1$, and hence $st(f)-1 \leq reg(f)$. 

To prove the other inequality, let $I$ be the saturation of the Jacobian ideal $J_f$, as discussed in a more general setting in \cite{HS}. Consider the exact sequence of graded $S$-modules
$$ 0 \to I/J_f \to M(f) \to S/I \to 0.$$
Then \cite[Corollary 20.19]{EisCA} implies that
$$reg(f) =reg M(f) \leq \max\{reg(I/J_f),reg(S/I)\}.$$
Note  the module $I/J_f$ has finite length, so \cite[Corollary 4.4]{Eis} implies that
$$reg(I/J_f)=\max \{k : (I/J_f)_k \ne 0\}=sat(J_f)-1,$$
in the notation from \cite{DBull}. Moreover, \cite[Corollary 2]{DBull} says that
$$sat(J_f) \leq \max \{T-ct(f),st(f)\},$$
where $T=3(d-2)$. On the other hand the quotient $S/I$ is a Cohen-Macaulay module satisfying $\depth S/I= \dim S/I =1$, 
and \cite[Corollary 4.8]{Eis} tells us that
$s=reg(S/I)$, where $s$ is the smallest integer such that $k \geq s$ implies
$$H_{S/I}(k)= P_{S/I}(k)= \tau(C).$$
This integer $s$ is determined in \cite[Proposition 2]{DBull}, where it is shown that one has $s=T-ct(f)$.
It follows that 
$$reg(f) \leq \max \{ T-ct(f), st(f)-1\}.$$
When $C:f=0$ is free, we have $ct(f)+st(f)=T$ by \cite{Dmax}, and hence we get $reg(f) \leq st(f).$
A direct computation using the definition of the Castelnuovo-Mumford regularity in terms of a resolution, see \cite{EisCA}, p. 505, yields  $reg(f)=d_2+d-3$ when $C$ is free.
The equality $st(f)=reg(f)$ follows using the above formulas for $st(f)$. Note that in the free case, the equality 
$reg(f)=st(f)$ is also a consequence of Theorem 4.2 in \cite{Eis}, since $M(f)$ is Cohen-Macaulay in this case.

When $C:f=0$ is not free, then it is shown in \cite[Corollary 1.7]{Dmax} that $ct(f)+st(f) \geq T+2$, which implies that $\max \{ T-ct(f), st(f)-1\}=st(f)-1$, and this completes the proof.

\endproof

\begin{cor}
\label{corG} 
Let $C:f=0$ be a reduced plane curve of degree $d\geq 4$. Then the following hold.

\begin{enumerate}

\item $C$ is free if and only if 
$reg(f)=2(d-2)-mdr_e(f)$. 

\item $C$ is nearly free if and only if 
$reg(f)=2(d-2)-mdr_e(f)+1$.

\item If $C$ is neither free nor nearly free, then
$reg(f) \geq 2(d-2) - mdr_e(f)+2$.

\end{enumerate}
\end{cor}

\proof  The claims (1) and (2) then follow from the equalities $ct(f)+st(f)=T$ (resp. $ct(f)+st(f)=T+2$) which are shown in  \cite[Corollary 1.7]{Dmax} to characterize the free (resp. nearly free) curves.
The claim in (3) follows from the equality \eqref{REL_e} and the inequality $ct(f)+st(f) \geq T+3$
which holds in this case by \cite[Corollary 1.7]{Dmax}. 
\endproof
Note that in the cases (1) and (2) above one has $mdr(f)=mdr_e(f)$, while in the case (3) both
$mdr(f)=mdr_e(f)$ and $mdr(f)=d-1 < mdr_e(f)$ may occur.

\begin{cor}
\label{corB} If $\A:f=0$ is any  arrangement of $d\geq 4$ lines, then  
$$reg(f) \leq 2d-5 \text{ and }
st(f)\leq 2d-4.$$
When $\A:f=0$ is a generic  arrangement of $d\geq 4$ lines, both of the above inequalities become equalities.
\end{cor}
\proof
The inequality $reg(f) \leq 2d-5$, with equality for a generic arrangement, follows from \cite[Corollary 3.5]{SchCMH}. The reader must notice that the regularity there is the Castelnuovo-Mumford regularity of the graded $S$-module $AR(f)$, and not as in our paper that of the Milnor algebra $M(f)$. The exact sequence \eqref{EXSE}
allows us to pass from one regularity to the other, namely one has
$$reg(f)=reg(M(f))=reg(AR(f))+d-3.$$
The formula for $st(f)$ in the case of a generic arrangement follows from 
\cite[Corollary 1.3]{DStEdin}. In fact, using Theorem \ref{thmA}, we need only one of these two invariants, since it is known that such a line arrangement is not free, see for instance
\cite{DS14} or Proposition \ref{proptau} below.

In the general case, if the arrangement $\A$ is not free, then $st(f)=reg(f)+1 \leq 2d-4.$ And for a free arrangement $\A:f=0$, one has $st(f)=d+d_2-3 \leq 2d-4$ since clearly $d_2 \leq d-1$.
\endproof

For the following result we refer to \cite{Dcurvearr}. For the case of hyperplane arrangements in $\PP^n$ see \cite{Yu}.
 \begin{thm}
\label{thmF} 
The set $F(d,\tau)$ of free curves in the variety $\CC(d,\tau) \subset \CC(d)_0$ of reduced plane curves of degree $d$
with a fixed global Tjurina number $\tau$ is a Zariski open subset in $\CC(d,\tau)$.

The set $FA(d,\tau)$ of free line arrangements in the variety $A(d,\tau) \subset A(d)$ of line arrangements consisting of $d$ lines and 
with a fixed global Tjurina number $\tau$ is a Zariski open subset in $A(d,\tau)$.
\end{thm}

Using the Galois covering $\psi:X(d) \to A(d)$, we introduce the following notation.

\begin{equation}
\label{not1}
X(d, m_k \leq m)=\psi^{-1}(\{ f \in A(d)  :  m_k(f) \leq m\} )  
\end{equation}
and $X(d, m_k = m)=X(d, m_k \leq m) \setminus X(d, m_k \leq m-1).$
\begin{equation}
\label{not2}
X(d, \tau \leq m)=\psi^{-1}(\{ f \in A(d) :  \tau(f) \leq m\} ) \text{ and } X(d,\tau)=\psi^{-1}(A(d,\tau)).
\end{equation}
\begin{equation}
\label{not3}
X(d, mdr \leq m)=\psi^{-1}(\{ f \in A(d)  : mdr(f) \leq m\} )
\end{equation}
and $X(d, mdr = m)=X(d, mdr \leq m) \setminus X(d, mdr \leq m-1).$
\begin{equation}
\label{not4}
FX(d, \tau)=\psi^{-1}(FA(d,\tau) )  \text{ and } FX(d)= \cup_{\tau}FX(d,\tau) .
\end{equation}

Then we have the following obvious consequence of Proposition \ref{prop2} and Theorem \ref{thmF}.
\begin{cor}
\label{cor1} The sets $X(d, m_k \leq m)$  and $X(d, \tau \leq m)$ are Zariski open in $X(d)$ for any
positive integers $k$ and $m$. The set $X(d; mdr \leq m)$ is Zariski closed in $X(d)$ for any positive integer $m$. Moreover, the set $FX(d,\tau)$ is Zariski open in the variety $X(d,\tau)$.
\end{cor}

It is known that $FA(d, \tau) \ne \emptyset$ implies that there is an integer $r \in [0,(d-1)/2]$ such that 
\begin{equation}
\label{tau}
FX(d, \tau) \subset X(d, mdr =r), \text{ where }\tau =\tau(d,r)= (d-1)^2-r(d-1-r),
\end{equation}
 see \cite{Dmax}, \cite{DStFD}. This implies the following result.
\begin{cor}
\label{taumincor}
Let $\A: f=0$ be a free line arrangement in $A(d, \tau)$ and define
$$ \tau(d)_{min}=\frac{3}{4}(d-1)^2$$
for $d$ odd, and
$$\tau(d)_{min}=\lfloor \frac{3}{4}(d-1)^2 \rfloor +1$$
for $d$ even. Then $\tau= \tau(f) \geq \tau(d)_{min}$.

\end{cor}

\proof
The inequalities follow from the formula for $\tau(d,r)$ given above. Here $\lfloor x \rfloor $ denotes the integral part of the real number $x$.
\endproof

\begin{rk}
\label{rknf} It is shown in \cite{Dmax} that a line arrangement $\A: f=0$ with $d=|\A|$ is nearly free with exponents $d_1 \leq d_2=d-d_1$ if and only if 
$$\tau(\A)= \tau(d,d_1)-1.$$
\end{rk}

\subsection{Three group actions on parameter spaces}

Consider the connected algebraic group $G=Aut(\PP^2)=PGL(3,\C)$ of dimension 8.
This group acts naturally on the variety $\CC(d)$ and all the subsets $\CC(d)_0$, $A(d)$,
$F(d,\tau)$, $FA(d,\tau)$ are in fact $G$-invariant, hence they inherit a natural $G$-action, and are unions of $G$-orbits $G \cdot f$, for some $f \in \CC(d)$. Moreover, $G$ acts also on the variety $X(d)$ in a diagonal way, and such that the map $\psi: X(d) \to A(d)$ is $G$-equivariant. It follows that all the subsets 
$X(d, m_k \leq m)$, $X(d, \tau \leq m)$, $X(d,\tau)$,  $X(d; mdr \leq m)$, $FX(d, \tau)$ are also $G$-invariant, so they consists of unions of $G$-orbits, denoted by $G \cdot (\ell_1,...,\ell_d)$, for some $(\ell_1,...,\ell_d) \in X(d)$.

 \begin{prop}
\label{prop3} Let $(\ell_1,...,\ell_d) \in X(d)$ and denote $f=\psi(\ell_1,...,\ell_d) \in A(d)$. Then one has the following.

\begin{enumerate}

\item $\dim G \cdot (\ell_1,...,\ell_d)= \dim G \cdot f$.

\item $\dim G \cdot f=2$ if $d=1$, $\dim G \cdot f=4$ if $d=2$ and for $d=3$, one has $\dim G \cdot f=5$ if $\A:f=0$ consists of 3 concurrent lines, and $\dim G \cdot f=6$ if $\A:f=0$ consists of a triangle.

\item For $d\geq 4$, $\dim G \cdot f=5$ if $mdr(f)=0$, $\dim G \cdot f=7$ if $mdr(f)=1$
and $\dim G \cdot f=8$ if $mdr(f)>1$.

\end{enumerate}

\end{prop}

\proof 

The first claim follows since the map $\psi$ has finite fibers. To prove (2) and (3), note that one has 
$$\dim G \cdot f= \dim G -\dim Fix(f),$$
 where $Fix(f)$ is the stabilizer subgroup of $f$. The Lie algebra of $Fix(f)$ is exactly $AR(f)_1$, i.e. the linear Jacobian syzygies, see \cite[Proposition 1.1]{dPCTC}.
When $d=1$, we can take $f=x$ and it follows that $\dim AR(f)_1=6$, since in the notation from \eqref{rel_m} one takes $a=0$ and $b,c \in S_1$ arbitrary.
When $d=2$, we can take $f=xy$, and it follows that $\dim AR(f)_1=4$, since $a=\lambda x$, $b=-\lambda y$ and $c \in S_1$ arbitrary.
When $d=3$ there are two possibilities. The first one is $f=x^3+y^3$, when $\dim AR(f)_1=3$, as $a=b=0$ and $c \in S_1$. The second case is $f=xyz$ and then $\dim AR(f)_1=2$, since $AR(f)_1$ is spanned in this case by $(x,-y,0)$ and $(x,0,-z)$.

Assume now that $d \geq 4$. Then, if $r=mdr(f)=0$, this means that $\dim AR(f)_0=1$, which implies $\dim AR(f)_1=3$. If $r=1$, then it follows from \cite[Proposition 2.2]{dPCTC} that
$\dim AR(f)_1=1$. When $r>1$, one has $AR(f)_1=0$, so the claims in (3) are now proved.

\endproof

Let $\G$ be the Galois group of $\C$ over $\Q$. Then $\G$ acts on the parameter spaces $A(d)$ and $X(d)$ by acting on the coefficients of the defining equations. It follows that all the subsets 
$X(d, m_k \leq m)$, $X(d, \tau \leq m)$, $X(d,\tau)$,  $X(d; mdr \leq m)$, $FX(d, \tau)$, as well as $A(d, m_k \leq m)$, $A(d, \tau \leq m)$, $A(d,\tau)$,  $A(d; mdr \leq m)$, $FA(d, \tau)$
are also $\G$-invariant.

The symmetric group $\sigma_d$ also acts on $X(d)$ by permuting the linear factors of the defining equation $f=0$ of a line arrangement, and this is the reason why some strata in $X(d)$ are not irreducible, while their images in $A(d)$ have this property, see for instance Proposition \ref{prop5} (1).

\subsection{On rigid plane curves and line arrangements}

We say that a plane curve $C:f=0$ is algebraically rigid if $(I/J_f)_d=0$,
where $I$ denotes as above the saturation of the Jacobian ideal $J_f$.
Indeed, the vector space $(I/J_f)_d$ is naturally identified to the space of first order locally trivial deformation of $C$  in $\PP^2$, modulo the above $G$-action, see \cite{Se1}, \cite{Se2}. These deformations preserve the analytic isomorphism type of each singular point of $C$.

\begin{ex}
\label{exalgrigid} 
It is known that a curve $C:f=0$ is free if and only if $I=J_f$, see \cite{ST}. In particular, any free curve is algebraically rigid. A generic line arrangement of 4 lines in $\PP^2$ is not free, but it is algebraically rigid
by Proposition \ref{prop3} (3) since $mdr(f)>1$ in this case. In fact, in this case one has $\dim (I/J_f)_3=1$ and $(I/J_f)_k=0$ for $k \ne 3$.
\end{ex}

We say that a reduced plane curve $C:f=0$ is topologically rigid if any deformation of $C$ preserving the number of irreducible components of $C$, their degrees and the topological type of each singularity of $C$ is trivial modulo the above $G$-action. For more on this type of rigidity see
\cite{KuSh}.

\begin{ex}
\label{extoprigid} 
A line arrangement $\A$ consisting in $d \geq 4$ lines passing through one point satisfies $mdr(f)=0$, and it is free. Hence $\A$ is an algebraically rigid curve. On the other hand, we can modify the cross-ratio of a subset of 4 lines in $\A$ by moving one line, without changing the topology of the singularity, and hence such a family will not be contained in one $G-$orbit. Hence $\A$ is not topologically rigid.
\end{ex}

\begin{rk}
\label{rkrigidity} 
If a reduced plane curve has only simple singularities of type $A_k$, $D_k$ and $E_6$, $E_7$, $E_8$, then $C$ is algebraically rigid if and only if $C$ is topologically rigid. Indeed, for a simple singularity, a topologically constant deformation is the same as an analytically constant deformation. In particular, for a line arrangement $\A$ having only double and triple points, the two rigidity notions coincide. In such a case we will simply say that $\A$ is rigid. For examples of this situation, see Remark \ref{rkrigidZiegler} and the stratum $A(L(\Delta))$ in Proposition \ref{propd=6} below. 

\end{rk}
For more on the interest of rigidity in the study of line arrangements, see \cite{Abe+}.

\section{A partition of the parameter space $X(d)$}

From now on in this paper we assume that $d \geq 4$. For a fixed integer $d \geq 4$, we denote by $\LL(d)$ the set of all possible intersection lattices $L(\A)$, for line arrangements in $\PP^2$ consisting of $d$ distinct lines. 
For a lattice $L \in \LL(d)$, we denote by $X(L)$ the set of all elements $(\ell_1,...,\ell_d) \in X(d)$ such that the line arrangement
\begin{equation}
\label{XL1}
\A: \ell_1= ...= \ell_d=0
\end{equation}
has an intersection lattice $L(\A)$ isomorphic to $L$, see \cite{OT} for more on intersection lattices. We also set $A(L)=\psi(X(L)).$
Such a lattice gives in particular information on the multiple points $p$ in the arrangement $\A$, and about their multiplicities, denoted by $m_p \geq 2$. In particular, we define 
\begin{equation}
\label{XL1.5}
\tau(L)=\sum_p (m_p-1)^2= \tau(\A).
\end{equation}
By definition, we have the following partitions
\begin{equation}
\label{XL2}
X(d)= \cup_{L \in \LL(d)}X(L) \text{ and } X(d,\tau)= \cup_{L \in \LL(d), \tau(L)=\tau}X(L)
\end{equation}
and similarly
\begin{equation}
\label{AL2}
A(d)= \cup_{L \in \LL(d)}A(L) \text{ and } A(d,\tau)= \cup_{L \in \LL(d), \tau(L)=\tau}A(L).
\end{equation}

One has the following.
\begin{prop}
\label{prop4} For any lattice $L \in \LL(d)$, the following hold.

\begin{enumerate}

\item The sets $X(L)$ and $A(L)$ are constructible; they are also $G$-invariant and $\G$-invariant.

\item $X(L) \subset X(d,\tau(L))$ and $A(L) \subset A(d,\tau(L))$.

\item The function $mdr : A(L) \to \N$ attains its minimal value on a Zariski closed subset $F$ of $A(L)$, and in general $F \ne A(L)$.

\item The function $m_k : A(L) \to \N$ attains its minimal value on a Zariski open subset $U_k$ of $A(L)$, and in general $U_k \ne A(L)$.

\end{enumerate}

\end{prop}

\proof The claim about the constructibility in (1) can be settled as follows. A point $p$ of multiplicity $k \geq 3$ will give rise to a set $\E(L)_p$ of $k-2$ equations to be satisfied by the set of coefficients $(a_i,b_i,c_i) \in \PP(S_1)$, where $\ell_i=a_ix+b_iy+c_iz$. Indeed, if the lines passing through $p$ are for instance $L_i: \ell_i=0$ for $i=1,2,...,k$, then the fact that all these lines pass through $p$ is expressed by the vanishing of $k-2$ determinants $D(1,2,j)$ of $3 \times 3$ matrices $A(1,2,j)$, constructed using the coefficients of $\ell_1, \ell_2$ and $\ell_j$ to define the corresponding three rows,
where $j=3,4,...,k$. Note that such determinants really define hypersurfaces in the product $\PP(S_1)^d$.
Moreover, when three lines $L_u$, $L_v$ and $L_w$ are not concurrent, we should add the condition that the corresponding determinant $D(u,v,w)$ is not zero. 
More details on this construction can be found in \cite{NY}, see however Remark \ref{rk0} below.
The $G$-invariance and the $\G$-invariance of $X(L)$ is obvious.

The claim (2) is clear. For the first part in claim (3), 
use Proposition \ref{prop2}, (2). 
For the second part of claim (3),
one may consider the example of two line arrangements 
$$\A: f=xy(x-y-z)(x-y+z)(2x+y-2z)(x+3y-3z)(3x+2y+3z)$$
$$(x+5y+5z)(7x-4y-z)=0$$
and
$$\A':f'=xy(x+y-z)(5x+2y-10z)(3x+2y-6z)(x-3y+15z)$$
$$(2x-y+10z)(6x+5y+30z)(3x-4y-24z)=0,$$
having isomorphic intersection lattices and constructed by Ziegler in \cite{Zi}. A picture of these arrangements can be found in \cite[Chapter 8]{DHA}.
They consists both of nine lines, and have only double and triple points. More precisely, they have $n_2= 18$ double points and $n_3=6$ triple points, and hence $\tau(\A)=\tau(\A')=42$. In the case of $\A$, the six triple points are on a conic, and a direct computation shows that 
$mdr(f)=5.$ 
For $\A'$,  the six triple points are not on a conic, i.e. the arrangement $\A'$ is a small deformation of the arrangement $\A$, and a direct computation shows that 
$mdr(f)=6$.
See also \cite[Example 13]{Schenck}.
The above example settles also the claim (4) by taking $k=13$, since
$$m_{13}(f')= \tau(f')=42 < m_{13}(f).$$
\endproof

\begin{rk}
\label{rkrigidZiegler} 
In fact, it is clear that there is a topologically constant 1-parameter family of line arrangements
$\A_t$ such that $\A_0=\A$ and $\A_t$ for $t \ne 0$ has the same numerical invariants as $\A'$. This family is obtained by moving the sixth triple point till it gets onto the conic determined by the first 5 triple points.
It follows that $\A$ is not rigid, and one can check that $\dim (I/J_f)_9=4$. A direct computation shows that
for $\A'$ one has $\dim (I/J_f)_9=4$ as well, i.e. $\A'$ is not rigid either.
\end{rk}

\begin{rk}
\label{rk-1} 

The above result says that the invariants $mdr$ and $m_k$ are not determined by the combinatorics in general. However, if $\A$ is a free arrangement, both $mdr(f)$ and $m_k(f)$ are determined by the lattice $L(\A)$. The claim for $mdr(f)$ follows from the formula
\eqref{tau}. The claim for $m_k(f)$ follows from the fact that the exponents $d_1=mdr(f)$ and $d_2=d-1-mdr(f)$ determine the Hilbert function $H_{M(f)}$ via the resolution given in Definition \ref{deffree} (2).

\end{rk}

\begin{rk}
\label{rk0} 
The set of equations $\E(L)=\cup_p\E(L)_p$ defined above is smaller than the set of equations 
constructed in \cite{NY}, and which we call $\E'(L)$ here. Indeed, any point $p$ of multiplicity $k\geq 3$ contributes
$k-2$ equations to our set $\E(L)$, and ${k \choose 3}$ equations to the set  $\E'(L)$.
The two ideals $I(\E(L))$ and $I(\E'(L))$ are distinct. Indeed, the equations in $\E'(L)$ are linearly independent degree 3-forms, as each of them involves monomials in distinct set of variables. For instance the monomial $a_1b_2c_3$ occurs only in the equation associated to the triple of lines $(L_1,L_2,L_3)$, supposed to pass through a multiple point $p$.

On the other hand, it is clear that the two ideals $I(\E(L))$ and $I(\E'(L))$ both have $Y(L)={\overline X(L)}$ as zero set, and hence one has in particular
\begin{equation}
\label{codim1}
\codim X(L) = \codim Y(L) \leq \sum_p (m_p-2),
\end{equation}
where $\codim X(L)$ means the codimension of $X(L)$ in the corresponding $X(d)$.
For lattices $L$ coming from line arrangements with few lines, or of a reduced complexity, the above inequality is an equality, see for an example Proposition \ref{prop5} (1) below. However, the monomial arrangement
$$\A(m,m,3): f=(x^m-y^m)(x^m-z^m)(y^m-z^m)=0,$$
has $d=3m$, 3 points of multiplicity $m$ and $m^2$ points of multiplicity $3$.
It follows that 
$$\sum_p (m_p-2)=3(m-2)+m^2>6m=\dim X(3m)$$
for $m \geq 5$. Hence for these values of $m$, the inequality \eqref{codim1} is strict.
\end{rk}

\begin{rk}
\label{rk1} The variety $X(L)$ corresponds exactly to the variety of all ordered complex realizations $\Sigma^{ord}(\CC)$ of the ordered combinatorics $\CC^{ord}$
 considered in \cite{Artal06}, where $\CC^{ord}$ is the ordered combinatorial type associated to the lattice $L$ with a fixed numbering of the lines.
The quotient $X(L)/G$  is the ordered moduli space $\M^{ord}(\CC)$ considered in \cite{Artal06}. The variety $A(L)$ corresponds exactly to the variety of all  complex realizations $\Sigma(\CC)$ of the  combinatorics $\CC$
as  considered in \cite{Artal06}, while $A(L)/G$ is the moduli space $\M(\CC)$ of the combinatorics $\CC$. If $L$ is the lattice corresponding to the MacLane line arrangement, it follows from \cite[Example 1.7]{Artal06} that $X(L)$ is the union of two $G$-orbits and in particular is not connected, while $A(L)$ is just one $G$-orbit, and hence it is irreducible.

\end{rk}

\begin{ex}
\label{ex1}  For any $d$, we denote by $L_{gen}$ the lattice of the generic line arrangement
of $d$ lines. Then by the above description $X(L_{gen})$ is a Zariski open subset of $X(d)$, and hence $\dim X(L_{gen})= \dim X(d)=2d$. Moreover $\tau(L_{gen})={ d \choose 2}$ and 
in fact one has $X(d,\tau(L_{gen}))=X(L_{gen})$, i.e. any lattice $L \in \LL(d)$ with 
$\tau(L)=\tau(L_{gen})$ is in fact isomorphic to the lattice $L_{gen}$. To prove this claim, recall the formula
\begin{equation}
\label{eqtau1}
\sum_p {m_p \choose 2} = {d \choose 2},
\end{equation}
valid for any line arrangement, see for instance \cite{Hirz}. Since 
$${m_p \choose 2} \leq (m_p-1)^2$$
for any $m_p\geq 2$, and the equality holds if and only if $m_p=2$, the claim follows
using the formula \eqref{XL1.5}. This argument implies also that
$\tau(L) > \tau (L_{gen})$ for any  lattice $L \in \LL(d)$, $ L \ne L_{gen}$.

Moreover, in this case it follows that the function $mdr \circ \psi$ is constant on $X(L_{gen})$, and it takes the value $d-2$, see \cite[Theorem 4.1]{DStEdin}, as well as all the functions $m_k$, since one has $ct(f)=st(f)=2d-4$ in this case, recall Corollary \ref{corArr} (2) and Corollary \ref{corB}.
\end{ex}

This example can be generalized as follows. For any $m$ satisfying $2 \leq m \leq d$, let $L(d,m)$ denote the intersection lattice of a line arrangement in $A(d)$ having one point of multiplicity $m$ and only double points in rest. Note that $L(d,2)=L_{gen}$ for any integer $d$. 

\begin{prop}
\label{prop5} Assume that $d \geq 4$. Then the following hold.

\begin{enumerate}

\item The sets $X(L(d,m))$ and $A(L(d,m))$ are smooth of dimension $(2d-m+2)$. Moreover $A(L(d,m))$ is irreducible.

\item $\tau(L(d,m))=(m-1)^2+m+(m+1)+....+(d-1)={d \choose 2} +{m-1 \choose 2}.$

\item For $2m>d$, the function $mdr $ is constant on $A(L(d,m))$, and it takes the values $d-m$.

\item For $m=d$ and $m=d-1$, any arrangement in $A(L(d,m))$ is free. Any arrangement in $A(L(d,d-2))$ is nearly free.
For $2 \leq m \leq d-3$, any arrangement in $A(L(d,m))$ is neither free,  nor nearly free.

\item Any line arrangement $\A:f=0$ with $mdr(f)=0$ satisfies $L(\A)=L(d,d)$.
Any line arrangement $\A:f=0$ with $mdr(f)=1$ satisfies $L(\A)=L(d,d-1)$.

\item $A(L(d,d)) =A(d,\tau(L(d,d))$ and $A(L(d,d-1)) =A(d,\tau(L(d,d-1)))$.

\end{enumerate}

\end{prop}

\proof To get an arrangement in $A(L(d,m))$, we have first to fix a point $p \in \PP^2$, and then $m$ distinct lines passing through $p$. These choices are parametrized by $B=\PP^2 \times U$, where $U$ is an open subset in $(\PP^1)^m$. Note that $B$ is smooth of dimension $m+2$ and irreducible. The remaining $d-m$ lines are to be chosen in a Zariski open set $F \subset (\PP(S_1))^{d-m}$, which is  smooth of dimension $2d-2m$ and irreducible. In this way we have constructed a fibration $F \to A(L(d,m)) \to B$, proving the first claim (1). 

Note that $X(L(d,m))$ is not connected in general. Indeed, for $d=4$ and $m=3$, we cannot continuously deform within $X(L(4,3))$ an element 
$(\ell_1,\ell_2,\ell_3,\ell_4)$
where the lines $L_j:\ell_j=0$ are concurrent for $j=1,2,3$ to an element
$(\ell'_1,\ell'_2,\ell'_3,\ell'_4)$ where the lines $L'_j:\ell'_j=0$ are concurrent for $j=2,3,4$.

The second claim follows from the formula \eqref{XL1.5}. The third claim follows from \cite[Theorem 1.2]{Dcurvearr}. 
The  claim (4) follows from the formula for $\tau(L(d,m))$ given in (2). Indeed, if a line arrangement $\A:f=0$ in $A(d)$ is free, then one has
$\tau(\A)=\tau(d,r)$
where $r=mdr(f)$, as explained in \eqref{tau}.

Suppose first that $2m>d$, and hence $r=d-m$.
The formula for $\tau(L(d,m))$ given in (2) shows that 
\begin{equation}
\label{eqtau10}
\delta= \tau(d,r) -\tau(\A)=\frac{(d-m)(d-m-1)}{2}.
\end{equation}
Hence we have the equality $\delta=0$ only for $m=d$ or for $m=d-1$. 
Assume now that $r \leq d-m-1$. Then  \cite[Theorem 1.2]{Dcurvearr} implies that either $r=m-1$ or $m\leq r$.
In the first case the arrangement is free with exponents $d_1=m-1$ and $d_2=d-m$, and the equation
\eqref{eqtau10} shows that this is possible only if $m=d-1$ or $m=d$, which is impossible.
In the second case, \cite[Theorem 1.1]{Dcurvearr} shows that $\tau(\A) \leq \tau(d,m)$ and equality holds exactly when $r=m$ and $\A$ is free. A direct computation shows that
\begin{equation}
\label{eqtau11}
\delta'= \tau(d,m) -\tau(\A)=\frac{(d-m-2)(d-m-1)}{2}+(m-1)>0
\end{equation}
for $2 \leq m \leq d-2$. For claim involving the nearly free arrangements, use the above and Remark \ref{rknf}.

The first part in claim (5) is clear, since $mdr(f)=0$ if and only if $f$ does not depend on the variable $z$ after a coordinate change. If $mdr(f)=1$, then let $m \geq 2$ be the maximal multiplicity of an intersection point in $\A$. Using \cite[Theorem 1.2]{Dcurvearr}, we deduce that 3 cases are possible.

(a) The case $mdr(f) =d-m=1$, which clearly settles our claim.

(b) The case $mdr(f)=m-1$, impossible, since this would imply that $\A$ is a generic arrangement, for which $mdr(f)=d-2 \geq 2$.

(c) The case $m \leq mdr(f)$, which is clearly impossible.

For claim (6), let $\B:g=0$ be a line arrangement in $A(d,\tau(L(d,d))$ (resp. in $A(d,\tau(L(d,d-1))$). Then $mdr(g) \geq 0$ (resp. $mdr(g) \geq 1$) and the claim follows using
\cite[Theorem 1.1]{Dcurvearr}, which, though not stated there, holds for $mdr(g) \geq 0$ as well.

\endproof

\begin{cor}
\label{cor2} For $d\geq 4$ and $m \in \{2,d-1,d\}$, one has
$$A(L(d,m))=A(d, \tau(L(d,m))=A(d, mdr=d-m)$$
and
$$X(L(d,m))=X(d, \tau(L(d,m))=X(d, mdr=d-m).$$

\end{cor}

\bigskip
Here is another example of a simple combinatorics for line arrangements.
Let $\tilde{L}(m_1,m_2)$ be the lattice of a projective line arrangement $\A\; :\; f=0$ obtained by the generic intersection of two pencils of $m_1$, respectively $m_2$ lines, with $m_2\geq m_1\geq 2$. Hence, a corresponding arrangement has  $d=(m_1+m_2)$ lines, $m_1m_2$ double points, one point of multiplicity $m_1$ and one point of multiplicity $m_2$. 
We will prove that such an arrangement $\A$ is never free.

\begin{prop}
\label{propAD}  
With this notation, one has the following.
\begin{enumerate}

\item The set $A(\tilde{L}(m_1,m_2))$ is smooth, irreducible of dimension $d+4$.

\item $\tau(\tilde{L}(m_1,m_2))=(d-1)^2-m_1m_2+1.$

\item The function $mdr $ is constant on $A(\tilde{L}(m_1,m_2))$ and takes the value $m_1$.

\item The intersection $FA(d) \cap A(\tilde{L}(m_1,m_2))$ is empty.

\end{enumerate}

\end{prop}

\proof

The first claim can be proved by a similar argument as that used in the proof of Proposition 
\ref{prop5} (1). The second claim is obvious.

By \cite[Theorem 1.2]{Dcurvearr} (applied for $m=m_2$), either $mdr(f)=m_1$, or $mdr(f)\leq m_1-1$. In the second case, either $mdr(f)\leq m_2-1$ and $2m_2<m_1+m_2$, contradiction with $m_1\leq m_2$, or $m_2\leq mdr(f)\leq m_1-1$, again contradiction with $m_1\leq m_2$. In conclusion, 
$mdr(f)=m_1$ and this proves (3).
To prove (4), note that the formula for $\tau(\tilde{L}(m_1,m_2))$ can be rewritten in the form
$$\tau(\tilde{L}(m_1,m_2))=(m_1-1)^2+(m_2-1)^2+m_1m_2=(d-1)^2-m_1(m_2-1)+1-m_1.$$ Then, by \eqref{tau}, the arrangement $\A$ is not free, since $m_2-1=d-m_1-1$ and $m_1 >1$.

\endproof

\begin{ex}
\label{exrkAD}

A slight variation of the previous configuration yields examples of free arrangements. 
Let ${\hat{L}}(m_1,m_2)$ be the lattice of a projective line arrangement $\A\; :\; f=0$ having exactly one line containing one point of multiplicity $m_1$ and one point of multiplicity $m_2, \; m_2\geq m_1\geq 3$, only double points apart from that, and $d=|\A|=m_1+m_2-1$.

By \cite[Theorem 1.2]{Dcurvearr} (applied for $m=m_2$), either $mdr(f)=m_1-1$, or $mdr(f)\leq m_1-2$. In the second case, it can only happen that $mdr(f)\leq m_1-2\leq m_2-1$ and $2m_2<m_1+m_2$, contradiction with the assumption $m_1\leq m_2$. So, $mdr(f)=m_1-1$ and $\tau({\hat{L}}(m_1,m_2))=(m_1-1)^2+(m_2-1)^2+(m_1-1)(m_2-1)$. We already know that such an arrangement $\A$ is free, since it is supersolvable (see \cite[Prop 5.114]{OT} and \cite[Theorem 4.2]{JT}). See also \cite{DStexpo} for this family of line arrangements.

In general one has  $X({\hat{L}}(m_1,m_2)) \ne X(m_1+m_2,\tau({\hat{L}}(m_1,m_2)))$. Indeed, by  \cite{DStexpo} any exponents $2 \leq d_1 \leq d_2$ of a free line arrangement can be obtained by such an arrangement. But there are free arrangements $\B$ which are not of this type, e.g. the monomial arrangements $\A(m,m,3)$ for $m \geq 2$ considered in Remark \ref{rk0}.
\end{ex}

The classification of the line arrangements $\A:f=0$ with $mdr(f)=2$ is given by the following theorem,  which is one of the main  results in \cite{T}. We give a proof of this classification from a new viewpoint. 

\begin{thm}
\label{thmToh} Let $\A:f=0$ be a line arrangement in $\PP^2$, with $mdr(f)=2$. Then $d=|\A| \geq 4$ and $\A$ is one of the following type of line arrangements,
described by their intersection lattices.

\begin{enumerate}

\item $\A \in A(L(d,d-2))$, or

\item $\A \in A({\hat{L}}(3,d-2))$ with $d \geq 5$, or 

\item $\A$ is  linear equivalent to the monomial arrangement $\A(2,2,3)$.

\end{enumerate}

\end{thm}

\proof

 Let $m$ be the maximal multiplicity of an intersection point in $\A$. If we denote $d=|\A|\geq 2$ and we assume $mdr(f)=2$, then \cite[Theorem 1.2]{Dcurvearr} implies that only the following cases are possible.

\medskip

\noindent CASE 1. $mdr(f)=d-m$, in other words $m=d-2$. This case covers the two 
cases (i) and (ii) in \cite[Theorem 2]{T}. Indeed, the case (1) corresponds to the case when $\A$ has a triple point except the point of multiplicity $m=d-2$ (e.g. for $d=5$ we have two triple points in $\A$), while the case (ii) corresponds to the case when $\A$ has only double points  except the point of multiplicity $m=d-2$. In fact, the line arrangements of type (ii) are exactly the line arrangement in $A(L(d,d-2))$ considered in Proposition \ref{prop5} (4) above, in particular they are all nearly free.
When $d \geq 5$, the arrangements of type (i) are exactly the line arrangements with  the intersection lattice of type ${\hat{L}}(3,d-2)$ which is discussed in Example \ref{exrkAD} above, in particular they are all free.

\medskip

\noindent CASE 2. $mdr(f)=m-1$ and $\A$ is free. In particular, this implies that the exponents of $\A$ are $d_1=2 \leq d_2$, and hence $d\geq 5$. Moreover
 $m=3$, and hence $\A$ has $n_2$ double points, $n_3$ triple points and no points of multiplicity $>3$. Using the formulas \eqref{free1} and \eqref{eqtau1}, we get the equations
 $$n_2+4n_3=d^2-4d+7 \text{ and } 2n_2 + 6n_3=d^2-d.$$
 They imply that $n_2=10d-d^2-21 \geq 0$ which yields $d\leq 7$. It is easy to classify the free line arrangements with $5 \leq d \leq 7$ and only double and triple points and we get in this way the lattice ${\hat{L}}(3,3)$ already seen above, and
the case (3) in \cite[Theorem 2]{T}, which is essentially the monomial arrangement $\A(2,2,3)$.
 
 \medskip
 
 \noindent CASE 3. $m=mdr(f)$, in other words $\A$ has only double points, as in Example \ref{ex1} above. But one knows that in this case $mdr(f)=d-2$, see \cite[Theorem 4.1]{DStEdin}, and hence we get again the case (ii) from
\cite[Theorem 2]{T} for $d=4$. 

\endproof

The above Theorem and Proposition \ref{prop5} imply the following.
\begin{cor}
\label{cormdr=2}
Any line arrangement $\A:f=0$ with $mdr(f) \leq 2$ is either free or nearly free.
Moreover, when $mdr(f) \leq 2$, the lattice $L(\A)$ determines the values of $mdr(f)$ and whether $\A$ is free or nearly free. In fact one has the following, where $d=|\A|$.
\begin{enumerate}

\item $mdr(f)=0$ if and only if $L(\A)=L(d,d)$;

\item $mdr(f)=1$ if and only if $L(\A)=L(d,d-1)$;

\item $mdr(f)=2$ if and only if $L(\A)$ is one of the lattices $L(d,d-2)$, ${\hat{L}}(3,d-2)$ or $L(\A(2,2,3))$.

\end{enumerate}

\end{cor}

\begin{rk}
\label{rkmdr=3}
The line arrangements $\A:f=0$ with $mdr(f)=3$ can be classified using the same approach, but the number of possibilities is much higher. Moreover, there are line arrangements with $mdr(f)=3$ which are neither free, nor nearly free, for instance
the generic arrangement of 5 lines.
\end{rk}

\begin{ex}
\label{exrkexpo}
We introduce a final lattice type. 
 For two integers $i \leq j$ we define a homogeneous polynomial in $\C[u,v]$ of degree $j-i+1$ by the formula 
\begin{equation}
\label{gij}
 g_{i,j}(u,v)=(u-iv)(u-(i+1)v) \cdots (u-jv).
\end{equation}
Consider the line arrangement $\A : f=0$ of $d=m_1+m_2 \geq 4$ lines in $\PP^2$ given by
$$f(x,y,z)= x(y-z)g_{1,m_1-1}(x,y)g_{2,m_2}(x,z)=0$$
for $2 \leq m_1 \leq m_2$. Denote by $L'(m_1,m_2)$ the corresponding intersection lattice $L(\A)$. One can show that the following hold, see for instance \cite{DStexpo}.
\begin{enumerate}

\item The line arrangement $\A$ has one point of multiplicity $m_1$, one point of multiplicity $m_2$, in addition to $(m_1-2)$ points of multiplicity 3 and $m_1(m_2-3)+6$ nodes;

\item $mdr(f)=m_1$;

\item $\tau(f)= (d-1)^2-m_1(m_2-1)-1$.

\end{enumerate}

\end{ex}

\begin{rk}
\label{rkrigidLattice}
We say that a lattice $L$ is rigid if the corresponding constructible set $A(L)$ is the disjoint union of finitely many $G-$orbits. It is clear that if $\A:f=0$ corresponds to a point in $A(L)$ with $L$ rigid, any topologically constant deformation of $\A$ is in fact a path in the connected component of $A(L)$ containing $\A$, which is by definition a $G-$orbit.
It follows that any such line arrangement $\A$ is topologically rigid.
Notice that the lattice $L(d,m)$ is rigid for $d \geq 4$ if and only if either $(d,m)=(4,2)$ or $(d,m)=(4,3)$. This follows from Proposition \ref{prop3}
and Proposition \ref{prop5}.
Other examples of rigid lattices $L$ are given in the next section.

A case of special interest is when the Galois group $\G$ acts transitively on the set of orbits in $A(L)$ for a rigid lattice $L$, see \cite{Abe+},
\cite{10lines}, \cite{Artal06}.
\end{rk}

\section{On the partition $A(d)= \cup_{L \in \LL(d)}A(L)$} 

In this section we describe the partition $A(d)= \cup_{L \in \LL(d)}A(L)$ for $ 4 \leq d \leq 6$,
and show that the complexity of this partition increases rapidly with $d$.

\subsection{ The case $d=4$.} \label{ex2}

For  $d=4$, the list $\LL(4)$ consists of 3 lattices, namely $L(4,2)$, $L(4,3)$ and $L(4,4)$ in the notation from  Proposition \ref{prop5}. Hence we have the following partition
$$A(4)=A(L(4,2)) \cup A(L(4,3)) \cup A(L(4,4)),$$
where $\dim A(L(4,2))=8$, $\dim A(L(4,3))=7$ and $\dim A(L(4,4))=6$. Moreover, the sets
$A(L(4,2))$ and $A(L(4,3))$ are $G$-orbits, i.e. the corresponding arrangements are rigid, while $A(L(4,4))$ is the union of a 1-parameter family of $G$-orbits, as can be seen using Proposition \ref{prop3} and its proof. Recall also Example \ref{extoprigid}.
Note that the closure of $A(L(4,2))$ in $A(4)$ is the whole space $A(4)$, while the closure of $A(L(4,3))$ in $A(4)$ is $A(L(4,3)) \cup A(L(4,4))$, which follows from Corollary \ref{cor1}. Moreover the set $ A(L(4,4))$ is closed in $A(4)$.

In this case, one has 
$$6=\tau(L(4,2)) < 7=\tau(L(4,3)) < 9=\tau(L(4,4))$$
and hence the corresponding 3 strata are distinguished by their Tjurina numbers.
Moreover, one has
$$FA(4)=A(L(4,3)) \cup A(L(4,4)).$$
Note also that even in this simple case, the set
$$A(4, st  \leq 3):=\{ f \in A(4)  : st(f) \leq 3\} )=A(L(4,3))$$
is neither open nor closed. Hence the invariant $st$ does not have nice semicontinuity properties as $\tau$ or $mdr$. By inspection of this list, we can state the following result.

\begin{prop}
\label{propd=4}  
With this notation, one has the following complete list of free and nearly free line arrangements for $d=4$.
\begin{enumerate}

\item The set $A(4,6)$ is open, 8-dimensional, and consists only of nearly free arrangements with exponents
$d_1=d_2=2$.

\item The set $A(4, 7)$ is irreducible, 7-dimensional, and consists only of  free arrangements
with exponents
$d_1=1, d_2=2$.

\item The set $A(4, 9)$ is irreducible, 6-dimensional, and consists only of  free arrangements with exponents
$d_1=0, d_2=3$.

\end{enumerate}

\end{prop}

\subsection{ The case $d=5$.}
\label{ex3}  For  $d=5$, the list $\LL(5)$ consists of $L(5,2)$, $L(5,3)$, $L(5,4)$, $L(5,5)$ and an additional lattice $L=L(\A)$ where $\A: xyz(x+y)(x+z)=0$. Note that the lattice $L$  is just the lattice ${\hat{L}}(3,3)$ from Example \ref{exrkAD}.
In this case one has 
$$A(5,mdr=2)= A(L(5,3)) \cup A(L),$$
with $\dim A(L(5,3))=9$ and 
$\dim A(L)=8.$ Hence Corollary \ref{cor2} does not hold for $m=d-2$ in this case.
One also has 
$$10=\tau(L(5,2)) < 11=\tau(L(5,3)) < 12= \tau(L)< 13=\tau(L(5,4))< 16=\tau(L(5,5)). $$
Hence again the corresponding 5 strata are distinguished by their Tjurina numbers.
Moreover, one has
$$FA(5)=A(L) \cup A(L(5,4)) \cup A(L(5,5)).$$

Note that 
$$\overline {A(L)} \cap A(L(5,4))=\emptyset.$$
Though this might be obvious for some readers, we prefer to give an argument which is likely to work in many similar situation. Note that, using the Curve Selection Lemma, if 
$\overline {A(L)} \cap A(L(5,4)) \ne \emptyset$, then we get a deformation of an ordinary singular point $(Y_4,0)$ of multiplicity 4 into two ordinary singular points $(Y_3,0)$ of multiplicity 3. Such a deformation is impossible, since it would contradict the semicontinuity of the spectrum on the interval $I=(-1/3, 2/3]$, see for details \cite[Theorem (8.9.8)]{Ku}.
Indeed, one has
$$ 1=\deg_I spec (Y_4) < 2 \deg_I spec (Y_3)=2.$$
It follows that
$$\overline {A(L)} \cap \overline {A(L(5,4))}=A(L(5,5)).$$
 By inspection of the list of lattices in $\LL(5)$, we can state the following result.

\begin{prop}
\label{propd=5}  
With this notation, one has the following complete list of free and nearly free line arrangements for $d=5$.
\begin{enumerate}

\item The set $A(5,11)$ is  irreducible, 9-dimensional, and consists only of nearly free arrangements with exponents
$d_1=2, d_2=3$. Any $f \in X(5,11)$ satisfies $st(f)=6$ and $reg(f)=5.$

\item The set $A(5, 12)$ is irreducible, 8-dimensional, and consists only of  free arrangements with exponents
$d_1=d_2=2$. Any $f \in X(5,12)$ satisfies $st(f)=reg(f)=4.$

\item The set $A(5, 13)$ is irreducible, 8-dimensional, and consists only of  free arrangements with exponents
$d_1=1, d_2=3$. Any $f \in X(5,13)$ satisfies $st(f)=reg(f)=5.$

\item The set $A(5, 16)$ is irreducible, 7-dimensional, and consists only of  free arrangements with exponents
$d_1=0, d_2=4$. Any $f \in A(5,16)$ satisfies $st(f)=reg(f)=6.$

\end{enumerate}

\end{prop}
Note that among the above sets, only $A(5,12)$ is a $G-$orbit and hence the corresponding arrangements are rigid. The arrangements in 
$A(5,13)$ are algebraically rigid, but not topologically rigid, recall Proposition \ref{prop3}.
\begin{rk}
\label{rksemic}
It is clear that $A(5,12)$ is contained in the closure of $A(5,11)$, and that $A(5,16)$ is contained in the closure of $A(5,12)$. The values given above for $st(f)$ and $reg(f)$ show that these invariants do not enjoy simple semicontinuity properties as in Proposition \ref{prop2}.

\end{rk}

\subsection{ The case $d=6$.}
\label{ex3.2}
For  $d=6$, the list $\LL(6)$ consists of 10 lattices. We list them in increasing order of their Tjurina numbers. 

$\bullet$ For $\tau= 15$, we have only the lattice $L(6,2)$ as predicted by the general theory, recall Example \ref{ex1}. Moreover $X(L(6,2))$ is an open subset in the 12-dimensional smooth variety $X(6)$

$\bullet$ For $\tau= 16$, we have only the lattice $L(6,3)$ and the corresponding set    $X(L(6,3))$ has codimension 1 in $X(6)$.

$\bullet$ For $\tau= 17$, we have two lattices, namely $\tilde L(3,3)$ and a new lattice, say
$\tilde L'(3,3)$. These two lattices have each 2 triple points and 9 nodes, and the invariant $mdr$ takes the value 3 in both cases. In the lattice $\tilde L(3,3)$ the 2 triple points are not on a line of the corresponding arrangement, while in the lattice $\tilde L'(3,3)$ the 2 triple points are  on such a line. In conclusion the corresponding two sets $X(\tilde L(3,3))$
and $X(\tilde L'(3,3))$ are not distinguished by the numerical invariants considered in this paper.  Indeed, since $ct(f)=7$ and $st(f)=8$ in both cases, the invariants $m_k$'s also coincide for any $k$.

Both sets $X(\tilde L(3,3))$ and $X(\tilde L'(3,3))$ have codimension 2 in $X(6)$

$\bullet$ For $\tau= 18$, we have again two lattices, namely the lattice $L(6,4)$, having a point of multiplicity 4, 9 nodes and $mdr=2$ and the lattice $L'(3,3)$ introduced in Example \ref{exrkexpo}, and having 3 triple points, 6 nodes and $mdr=3$.

The set $X(L(6,4))$ has codimension 2 in $X(6)$, while the set $X(L'(3,3))$ has codimension 3 in $X(d)$.

Comparing the values of $\tau(L) \leq 18$ and the corresponding values of the invariant $mdr$, we conclude that there are no free arrangements in this range.

$\bullet$ For $\tau= 19$, we have again two lattices, namely the lattice $ {\hat L}(3,4)$, having one point of multiplicity 4, one triple point and 6 nodes, and the lattice $L(\Delta)$ corresponding to the arrangement 
$$\A:f=(x^2-y^2)(x^2-z^2)(y^2-z^2)=0,$$
and hence having 4 triple points and 3 nodes.
Both of the corresponding sets $X(L)$  contain only free arrangements with
 $mdr=2$. 

$\bullet$ For $\tau= 21$, we have only the lattice $L(6,5)$.

$\bullet$ For $\tau= 25$, we have only the lattice $L(6,6)$.
The properties of the last two lattices are discussed in Proposition \ref{prop5} (4),  (5) and (6). In particular, the sets $X(d,\tau)$ for $\tau= 19, 21, 25$ consist only of free arrangements, i.e. the last claim in Corollary \ref{cor1} holds in a stronger version.

As a conclusion, we can state the following result.

\begin{prop}
\label{propd=6}  
With this notation, one has the following complete list of free and nearly free line arrangements for $d=6$.
\begin{enumerate}

\item The set $A(6,18)$  has two irreducible components, namely $A(L(6,4))$ of dimension 10, and  $A(L'(3,3))$ of dimension 9; they consist  only of nearly free arrangements with exponents
$d_1=2, d_2=4$, and respectively $d_1=d_2=3$.

\item The set $A(6, 19)$ has two irreducible components, namely $A( {\hat L}(3,4))$ of dimension 9, and  $A(L(\Delta))$ of dimension 8;
they consist only of  free arrangements with exponents
$d_1=2, d_2=3$.

\item The set $A(6, 21)$ is irreducible, 9-dimensional, and consists only of  free arrangements with exponents
$d_1=1, d_2=4$.

\item The set $A(6, 25)$ is irreducible, 8-dimensional, and consists only of  free arrangements with exponents
$d_1=0, d_2=5$.

\end{enumerate}

\end{prop}

\medskip

Note that $Y=A(L(\Delta))$ is the only $G$-orbit in the list above, and hence 
consists only of rigid arrangements. Moreover $Z=\overline Y \setminus Y$ is a closed $G$-invariant subset of dimension $<8$. Using the same type of argument as in the case $d=5$ above, one can show that $Z \subset A(L(6,6)$.
Since $A(L(6,6))$  is 8-dimensional by Proposition \ref{prop5} (1), it follows that $Z$ is not a union of strata in the partition. This shows in particular that this partition is not Whitney regular, see \cite[Chapter 1]{D1} for basic facts on regular stratifications.

\begin{rk}
\label{rkd=7}
It is possible to extend this discussion to $d=7$. For $d<7$ we have seen that the closed set $X(d, \tau \geq \tau(d)_{min})$
contains only free arrangements. For $d=7$, one has $ \tau(d)_{min}=27$.
The new aspect occurring in this case is that the set
$X(7,\tau \geq 27)$ contain free arrangements  and one nearly free arrangement type with exponents $d_1=2, d_2=5$.

\end{rk}
The above leads us to ask the following.
\begin{question}
\label{qd=7}
Is it true that, for any $d \geq 4$, the closed set $X(d, \tau \geq \tau(d)_{min})$ in $X(d)$ contains only  nearly free and free arrangements with $d$ lines?
\end{question}
Unfortunately the answer to this question is negative. To see this, it is enough to consider the arrangements in $X(L(d,d-3))$ for $d\geq 11$. Then using Proposition \ref{prop5} it is easy to check that $\tau(L(d,d-3)) > \tau(d)_{min}$. On the other hand, the formula \eqref{eqtau10} and the characterization of free arrangements (resp. nearly free arrangements) by the property $\delta=0$ (resp. $\delta =1$) given in \cite{Dmax} show that any arrangement in $X(L(d,d-3))$ for $d\geq 11$ is neither free nor nearly free.

\section{On Terao's conjecture} 
With the above notation, this conjecture in the case of line arrangements can be stated as follows.
\begin{conj}
\label{conjT}  [Terao's Conjecture for the line arrangement $\A$]
Let $\A$ be a free line arrangement with $d=|\A|$. Then 
$$A(L(\A)) \subset FA(d).$$
Equivalently, $X(L(\A)) \subset FX(d).$
\end{conj}
Assume that $\A$ is free with exponents $d_1 \leq d_2$. Then the following are known.

 \begin{thm}
\label{thmT1} 
Terao's conjecture holds for the line arrangement $\A$ if one has either $d=|\A| \leq 12$, or $d_1 \leq 5$.
\end{thm}
For the proofs of this result, see \cite{Abe1}, \cite{Abe+}, \cite{FV}. The next result was proved in \cite{Dcurvearr}, but we give below a new, clearer proof. 
Another rapid proof can be obtained by combining Theorem 2.7 and Lemma 2.10 in \cite{Abe+}. Let $m(\A)$ denote the maximal multiplicity of an intersection  point in $\A$.
 \begin{thm}
\label{thmT2} 
Terao's conjecture holds for the line arrangement $\A$ if $m(\A) \geq d_1$.
\end{thm}
\proof
We apply first \cite[Theorem 1.2]{Dcurvearr} to the free arrangement $\A$. It follows that 
\begin{equation}
\label{eqT1}
d_1 \in \{d-m, m-1,m\},
\end{equation}
where $m=m(\A)$. It follows that $m \in \{d-d_1, d_1, d_1+1\}$.
Let now $\B \in A(L(\A))$, given by $g=0$, and note that $m(\B)=m$. Now we apply 
\cite[Theorem 1.2]{Dcurvearr} to the line arrangement $\B$.

If we are in the case $mdr(g)=d-m(\B)=d-m \in \{d_1,d-d_1, d-d_1-1\}$, it follows that
$mdr(g) \geq mdr(f)=d_1$. Since $\tau(\B)= \tau(\A)$, the result \cite[Theorem 1.1]{Dcurvearr} implies that $\B$ is free with the same exponents $d_1 \leq d_2$.

If we are in the case $mdr(g)=m(\B)-1$, then the arrangement $\B$ is free, so there is nothing to prove. Finally, we have to consider the case $m=m(\B) \leq mdr(g)\leq d-m-1$.
For $m=d-d_1$ this implies $d-d_1 \leq d_1-1$, which is impossible.
For $m=d_1$ or $m=d-d_1-1=d_2$, we are again in the case $mdr(g) \geq mdr(f)=d_1$, and we conclude as above.
\endproof
As explained in \cite{Dcurvearr}, this Theorem implies the following.
\begin{cor}
\label{corT1} With the above notation, one has the following.

\begin{enumerate}

\item Terao's conjecture holds for the line arrangement $\A$ if $m(\A) \geq d/2$.

\item Terao's conjecture holds for the line arrangement $\A$ if $ d_1\leq \sqrt {2d+1}-1$.

\end{enumerate}

\end{cor}

We end with a result saying that a free arrangement cannot have too many singularities.
\begin{prop}
\label{proptau}  
The intersection $FX(d) \cap X(L)$ is empty if 
\begin{equation}
\label{eqtau0}
\tau(L)  < \frac{3}{4}(d-1)^2.
\end{equation}
In particular, the inequality \eqref{eqtau0} holds  if
$$\sum_p(m_p-1)>\frac{(d+3)(d-1)}{4}$$
where $p$ runs through the set of multiple points of the lattice $L$, and $m_p \geq 2$ denotes the multiplicity of $p$.
\end{prop}

\proof 
The first claim follows from Corollary \eqref{taumincor}.

Then the formula for $\tau(L)$ given in \eqref{XL1.5} and the equality \eqref{eqtau1} imply that
$$\tau(L)=2{d \choose 2}-\sum_p(m_p-1)<2{d \choose 2}-\frac{(d+3)(d-1)}{4}=\frac{3}{4}(d-1)^2.$$
\endproof

\begin{ex}
\label{extau}  Assume that the line arrangement $\A$ is not generic, but has a lot of nodes, namely it has  $N> \frac{(d+3)(d-1)}{4}-2$ nodes, besides some other multiple points. Then $\A$ is not free by the above result, since there is at least one point $p$ with $m_p \geq 3$. When $d=7$, this says that
an arrangement $\A$ having at least $14$ nodes satisfies
$$\tau(\A)  < \tau(7)_{min}=\frac{3}{4}(d-1)^2=27$$
and hence it is not free. A detailed classification of the line arrangements $\B$ with $|\B|=7$, shows that there is a nearly free arrangement $\B$ having  11 nodes and one point of multiplicity 5 such that
$\tau(\B) = 27.$
Moreover, for all arrangements $\B'$ having 12 nodes (and some other multiple points), one has $ \tau(\B') <27$. Hence our bound is two units apart from  the optimal one in this case.

\end{ex}

\end{document}